\documentclass[11pt]{article}
\usepackage{inputenc}
\usepackage{enumerate}
\usepackage{amssymb}
\usepackage{amsmath}
\usepackage{mathrsfs}
\usepackage{amsthm}
\usepackage{tikz-cd}
\usepackage{mathpazo}

\usepackage{extarrows}
\usepackage{color}
\usepackage{setspace}
\usepackage{hyperref}
\usepackage{graphicx, subfig}
\usepackage[square, comma, sort&compress, numbers]{natbib}
\usepackage{setspace}

\usepackage{titlesec}
\usepackage{titletoc}
\usepackage[titletoc]{appendix}

\theoremstyle{definition}
\newtheorem{Def}{Definition}[section]

\theoremstyle{plain}
\newtheorem{Thm}[Def]{Theorem}
\newtheorem{Lem}[Def]{Lemma}
\newtheorem{Pro}[Def]{Proposition}

\usepackage[left=2.6cm, right=2.6cm, top=3cm, bottom=3cm]{geometry}
\usepackage[all]{xy}
\usepackage{bm}

\def\IR{\mathbb R}\def\IC{\mathbb C}

\def\A{\mathcal A}\def\E{\mathcal E}\def\K{\mathcal K}\def\L{\mathcal L}\def\H{\mathcal H}

\def\supp{\textup{supp}}

\def\prop{\textup{Prop}}
\def\diam{\textup{diam}}

\def\ox{\otimes}

\def\ox{\otimes}
\def\EG{\underline{\E}\Gamma}
\def\Ga{\Gamma}
\def\car{\curvearrowright}

\author{Liang Guo\and Qin Wang}
\title{A coarse geometric approach to the generalised Green-Julg Theorem}
\begin{document}
\maketitle

\begin{abstract}
In this note, we provide a proof of the generalised Green-Julg theorem by using the language of twisted localization algebras introduced by G.~Yu. This proof is for those who have interests in coarse geometry but not so familiar with $KK$-theory or $E$-theory.

\end{abstract}

\section{Introduction}

The Dirac-dual-Dirac method introduced by G.~Kasparov \cite{Kas88} is one of the most powerful methods when dealing with the strong Novikov conjecture for a countable discrete group. This method has been developed in recent 30 years (see Section 2 in \cite{Tu1999} or Section 5 in \cite{KS2003}). In \cite{GHT2000}, E.~Guentner et al introduced the generalised Green-Julg theorem which is the key step to the Dirac-dual-Dirac method. In \cite{HK2001}, N.~Higson and G.~Kasparov show that the Baum-Connes conjecture with coefficient holds for all a-T-menable groups and the generalised Green-Julg theorem play an important role in the proof.

In \cite{Yu1997}, G.~Yu introduced a notion of the localization algebra for a proper metric space $X$ whose $K$-theory play a role of the $K$-homology of $X$. Later in \cite{Yu2000}, G.~Yu adapts and develops a construction of \cite{HK2001} to a coarse geometric version by using the twisted algebras. This method for the coarse Baum-Connes conjecture is always called the geometric Dirac-dual-Dirac method.

In this paper, we shall prove that the equivariant version of geometric Dirac-dual-Dirac method for cocompact action $\Ga\car X$ is actually identified with the original Dirac-dual-Dirac method. The main result of this paper is as follows:

\begin{Thm}\label{main result}
Let $\A$ be a proper $\Ga$-$C^*$-algebra. The evaluation map
$$ev_*:\lim\limits_{d\to\infty}K_*(C^*_L(P_d(\Ga),\A)^{\Ga})\to\lim\limits_{d\to\infty}K_*(C^*(P_d(\Ga),\A)^{\Ga})$$
is an isomorphism. We also have the following commuting diagram:
$$\begin{tikzcd}
\lim\limits_{d\to\infty}K_*(C^*_L(P_d(\Ga),\A)^{\Ga}) \arrow[r, "ev_*"] \arrow[d, "\cong "'] & \lim\limits_{d\to\infty}K_*(C^*(P_d(\Ga),\A)^{\Ga}) \arrow[d, "\cong"] \\
RK^{\Ga}_*(\EG,\A) \arrow[r, "\mu_{\A}"]                       & K_*(\A\rtimes_r\Ga)  ,                 
\end{tikzcd}$$
where $RK^{\Ga}_*(\EG,\A)=\lim_{Y\subseteq\EG,Y\text{ cocompact}}KK^{\Ga}_*(C_0(Y),\A)$. As a corollary, the assembly map
$$\mu_{\A}:RK_*^{\Ga}(\EG,\A)\to K_*(\A\rtimes_r\Ga)$$
is an isomorphism.
\end{Thm}

The paper is organised as follows: In section 2, we recall the twisted algebras for a cocompact action $\Ga\car X$ introduced by G.~Yu. We also show the diagram in Theorem \ref{main result} commutes. In section 3, we show the evaluation map is an isomorphism. As a result the assembly map is an isomorphism which completes the proof of Theorem \ref{main result}.

\section{Twisted localization algebra and $KK$-theory}

Let $\Ga$ be a countable discrete group. Equip $\Ga$ with a left invariant, proper, word length metric which is unique under coarse equivalence. Let $X$ be a proper metric space. An \emph{isometric group action} of $\Ga$ on $X$ is a group homomorphism $\Ga\to Isom(X)$, where $Isom(X)$ is the group of all isometries on $X$. A group action is always denoted by $\Ga\car X$. The action $\Ga\car X$ is said to be \emph{proper} if for any bounded $K\subseteq X$, there are only finite $\gamma\in\Ga$ such that $\gamma K\cap K\ne\emptyset$. This action can induce an action $\alpha:\Ga\to Aut(C_0(X))$ by using the formula
$$(\alpha_{\gamma}(f))(x)=f(\gamma^{-1}x).$$
for any $\gamma\in\Ga$, $x\in X$ and $f\in C_0(X)$. Throughout this section, we shall always assume $X$ admits an isometric proper cocompact action $\Ga\car X$.

A separable $\Ga$-$C^*$-algebra $\A$ is \emph{proper} if there exists a second countable, locally compact, proper $\Ga$-space $M$ and an equivariant $*$-homomorphism from $C_0(M)$ to the center of the muiltiplier algebra of $\A$ such that $C_0(M)\cdot\A$ is dense in $\A$. For any $a\in\A$, the support of $a$ is defined to be the complement of $x\in M$ for which there exists $f\in C_0(M)$ such that $a\cdot f=0$ and $f(x)\ne 0$. Equip $M$ with a proper metric $d_M$, from now on, we shall always view $M$ as a proper metric space.

We will not talk about the equivariant geometric module here. One is referred to Chapter 4 and Chapter 5 in \cite{HIT2020} for the standard coarse-geometric terminolgy. In this paper, we shall view an operator as an infinite-degree matrix. Let $\H_0$ be an infinite-dimensional Hilbert space, $\ell^2(\Ga)$ be the Hilbert space with the left-regular representation. Denote $\K=\K(\ell^2(\Ga)\ox\H_0)$ be the $C^*$-algebra of all compact operators on $\ell^2(\Ga)\ox\H_0$. Then $\K$ is made into a $\Ga$-$C^*$-algebra with the action given by the conjugation action on $\ell^2(\Ga)$.

\begin{Def}\label{twisted Roe}
Let $Z\subseteq X$ be a countable dense subset. The algebraic Roe algebra with coefficient in $\A$, denoted by $\IC[X,\A]^{\Ga}$, is the set of all bounded functions $T$ from $Z\times Z$ to $\K\ox\A$ such that\begin{itemize}
\item[(1)]there exists $L>0$ such that for each $x\in Z$,
$$\#\{y\in Z\mid T(x,y)\ne 0\}\leq L\quad\text{and}\quad\#\{y\in Z\mid T(y,x)\ne 0\}\leq L;$$
\item[(2)]there exists $r>0$ such that $T(x,y)=0$ for all $d(x,y)>r$. The infimum of such $r$ is called the propagation of $T$, denoted by $\prop(T)$;
\item[(3)]for any bounded $B\subseteq X$,
$$\#\{(x,y)\in B\times B\mid T(x,y)\ne 0\}<\infty;$$
\item[(4)]$T$ is $\Ga$-invariant, i.e., $T(x,y)=\gamma^{-1}\cdot (T(\gamma x,\gamma y))$ for any $\gamma\in\Ga$ and $x,y\in Z$.
\end{itemize}\end{Def}

Let $H=\ell^2(Z)\ox\ell^2(\Ga)\ox\H_0\ox\A$ be the Hilbert $A$-module. Then there is a faithful non-degenerate representation $\varphi:C_0(X)\to \L(H)$ given by
$$\varphi(f)\cdot(\delta_z\ox\delta_{\gamma}\ox v\ox a)=f(z)(\delta_z\ox\delta_{\gamma}\ox v\ox a).$$
The $\Ga$-action on $H$ is given by the left-regular action on $\ell^2(Z)$ and $\ell^2(\Ga)$. The production and $*$-operation of $\IC[X,\A]^{\Ga}$ is given by 
$$(TS)(x,y)=\sum_{z\in Z}T(x,z)S(z,y)$$
$$T^*(x,y)=(T(y,x))^*.$$
There exists a canonical representation of $\IC[X,\A]^{\Ga}$ on $H$ by matrix product as in \cite{Yu2000}, which fits the $\Ga$-action on $H$ clearly.

\begin{Def}
The equivariant Roe algebra with coefficient in $\A$, denoted by $C^*(X,\A)^{\Ga}$, is defined to be the completion of $\IC[X,\A]^{\Ga}$ under the norm given by its representation on $H$.
\end{Def}

We should mention that the algebra we define above is a dense subalgebra of $\IC[\Ga,X,A]$ in \cite{KY2012}. The triple $(C_0(X),\Ga,\varphi)$ is clearly an admissible covariant system. The original idea of this construction comes from the twisted Roe algebra defined by G.~Yu in \cite{Yu2000}. However it is slightly different with Yu's original definition since we do not have the support condition. The next lemma shows that the support condition is implied by the condition that the action $\Ga\car X$ is cocompact.

\begin{Lem}\label{support condition}
For each $T\in\IC[X,\A]^{\Ga}$ and $\varepsilon>0$, there exists $T'\in\IC[X,\A]^{\Ga}$, $R>0$ such that the diameter of $\supp(T'(x,y))$ is bounded by $R$ and $\|T-T'\|\leq\varepsilon$.
\end{Lem}

\begin{proof}
Since $\Ga\car X$ is cocompact and $X$ is a proper metric space, there exists a bounded subset $K_0\subseteq X$ such that $\Ga\cdot K_0=X$. Denote
$$K=\{x\in X\mid d(x,K)\leq\prop(T)\}.$$
By condition (3) in Definition \ref{twisted Roe}, there are only finite $(x,y)\in K\times K$ such that $T(x,y)\ne 0$. For any $\varepsilon'>0$ and $(x,y)\in K$, set $T'(x,y)\in \K\otimes\A$ such that $\supp(T'(x,y))$ is bounded and $\|T(x,y)-T'(x,y)\|<\varepsilon'$. Moreover, if $(x,y)$ and $(\gamma x,\gamma y)$ are both in $K\times K$, we assume that $T'(x,y)=\gamma^{-1}T'(\gamma x,\gamma y).$ We can always do that because $C_0(M)\A$ is dense in $\A$ and $T$ is $\Ga$-invariant. Take $R>0$ such that $\diam(\supp(T'(x,y)))\leq R$ for all $(x,y)\in K\times K$.

For any $(x,y)\notin K\times K$ with $d(x,y)\leq\prop(T)$, there must exists $\gamma\in\Ga$ such that $\gamma x\in K_0$. Then
$$d(\gamma y,K_0)\leq d(\gamma y,\gamma x)\leq \prop(T),$$
this means that $(\gamma x,\gamma y)\in K\times K$. Set $T'(x,y)=\gamma^{-1}T'(\gamma x,\gamma y)$. Set $T'(x,y)=0$ for all $(x,y)$ with $d(x,y)>\prop(T)$.

By \cite[Lemma 3.4]{GWY2008}, take $\varepsilon'$ small enough such that $\|T-T'\|\leq \varepsilon$. Then we complete the proof.
\end{proof}

The following proposition reveals the relationship of equivariant Roe algebra and the crossed product:

\begin{Pro}[\cite{KY2012}]\label{cp and Roe}
Denote $\A\rtimes_r\Ga$ the reduced crossed product algebra, then
$$K_*(C^*(X,\A)^{\Ga})\cong K_*(\A\rtimes_r\Ga).$$
\end{Pro}

Let $K_0$ be as in the proof of Lemma \ref{support condition}, i.e., $X=\Ga\cdot K_0$. Let $m_0\in M$ be a based point.

\begin{Def}
The localization algebra with coefficient in $\A$, denoted by $\IC_L[X,\A]^{\Ga}$, is the set of all bounded and uniform continuous functions $g$ from $\IR_+$ to $\IC[X,\A]^{\Ga}$ such that\begin{itemize}
\item[(1)] there exists $R>0$ such that $\diam(\supp(g(t)(x,y)))\subseteq B(m_0,R)$ for all $x\in K_0$;
\item[(2)] $\prop(g(t))$ tends to $0$ as $t$ tends to infinity.
\end{itemize}
The localization algebra, denoted by $C^*_L(X,\A)^{\Ga}$, is defined to be the completiong of $\IC_L[X,\A]^{\Ga}$ under the norm
$$\|g\|=\sup_{t\in\IR_+}\|g(t)\|.$$
\end{Def}

In \cite{KY2012}, G.~Kasparov and G.~Yu introduced a localization algebraic approach to equivariant $KK$-theory. Notice that there is a support condition in our definition of localization algebra which makes our localization algebra more like the twisted algebra defined in \cite{Yu2000}. Thus there is a canonical map
$$K_*(C^*_L(X,\A)^{\Ga})\to KK_*^{\Ga}(C_0(X),\A),$$
where the map is induced by the inclusion from our localization algebra to the localization algebra without the support condition. We shall show in the following proposition that the $K$-theory of our localization algebra also realise the equivariant $KK$-group.

\begin{Pro}
If $X$ is $\Ga$-compact, then $K_*(C^*_L(X,\A)^{\Ga})\cong KK^{\Ga}_*(C_0(X),\A)$.
\end{Pro}

\begin{proof}
First of all, we claim that
$$K_*(C^*_L(\Ga\times_{F}U,\A)^{\Ga})\cong K_*(C^*_L(U,\A)^F)$$
for each finite subgroup $F\leq\Ga$ and $F$-space $U$, where $\Ga\times_FU$ is the balanced product. Notice that $\Ga\times_FU=\bigsqcup_{[\gamma F]\in\Ga/F}\gamma U$. Let $r=\min\{d(U,\gamma U)\}$. For each $[g]\in K_0(C^*_L(\Ga\times_{F}U,\A)^{\Ga})$, since the propagation of $g$ tends to $0$ as $t$ tends to infinity, there exists $T>0$ such that $\prop(g(t))<\frac r2$. Define
$$g_U(t)=\chi_Ug(t+T)\chi_U\in K_0(C^*_L(U,\A)^F).$$
On the other hand, for any $h\in K_0(C^*_L(U,\A)^F)$, define
\[(h^{\Ga}(t))(x,y)=\left\{\begin{aligned}&0,&&x,y\text{ belongs to different component }\gamma U;\\
&\gamma^{-1}T(\gamma x,\gamma y),&&\gamma x,\gamma y\text{ both belong to }U\text{ for some }\gamma.\end{aligned}\right.\]
It is clear that the maps $[g]\mapsto [g_U]$ and $[h]\mapsto [h^{\Ga}]$ are group homomorphisms and they are the inverse to each other. The case for $K_1$ follows from a suspension argument. This shows the claim.

Both $KK$-theory and $K$-group of the localization algebra admit Mayer-Vietoris sequence and invariant under strong Lipschtiz homotopy equivalence (see \cite{Black} for $KK$-theory and \cite{Yu1997} for localization algebra). Since $X$ is cocompact, by using five lemma and a Mayer-Vietoris argument, it suffices to prove
$$K_*(C^*_L(\Ga\times_FU,\A)^{\Ga})\cong KK_*^{\Ga}(C_0(\Ga\times_FU),\A)$$
is an isomorphism for a finite subgroup $F\leq\Ga$ and a contractible subset $U\subseteq X$, the map is given by the inclusion. Fix an $m_0\in M$. Then $\Ga\cdot B(m_0,R)$ is a $\Ga$ invariant subspace of $M$. Denote
$$\A_R=\overline{C_0(\Ga\cdot B(m_0,R))\A}.$$
Let $C^*_L(\Ga\times_FU,\A_R)^{\Ga}$ be the subalgebra of $C^*_L(\Ga\times_FU,\A)^{\Ga}$ generated by all elements $g$ such that
$$\supp((g(t))(x,y))\subseteq B(m_0,R)$$
for all $x\in U$ and $t\in\IR_+$. For any $R'>R$, there is a canonical inclusion $C^*_L(\Ga\times_FU,\A_R)^{\Ga}\to C^*_L(\Ga\times_FU,\A_{R'})^{\Ga}$. It is clear that $\lim\limits_{R\to\infty}C^*_L(\Ga\times_FU,\A_R)^{\Ga}=C^*_L(\Ga\times_FU,\A)^{\Ga}$. Then it suffices to prove
$$\lim_{R\to\infty}K_*(C^*_L(\Ga\times_FU,\A_R)^{\Ga})\to KK^{\Ga}_*(C_0(\Ga\times_FU),\A)$$
is an isomorphism.

Notice that we have the following commuting diagram
$$\begin{tikzcd}
\lim\limits_{R\to\infty}K_*(C^*_L(\Ga\times_FU,\A_R)^{\Ga}) \arrow[r] \arrow[d, "\cong "'] & KK^{\Ga}_*(C_0(\Ga\times_FU),\A) \arrow[d, "\cong"] \\
\lim\limits_{R\to\infty}K_*(C^*_L(U,\A_R)^{F}) \arrow[r] \arrow[d, "\cong "']  & KK^{F}_*(C_0(U),\A)   \arrow[d, "\cong"]                 \\
\lim\limits_{R\to\infty}K_*(C^*_L(pt,\A_R)^F)\arrow[r]     &    KK^F_*(\IC,\A).
\end{tikzcd}$$
For each $R>0$, the group $K_*(C^*_L(pt,\A_R)^F)$ is isomorphic to $K_*(\A_R\rtimes_rF)$ by using the evaluation map (see \cite[Lemma 12.4.3]{HIT2020}) and Proposition \ref{cp and Roe}. Moreover, $KK_*^F(\IC,\A)\cong K_*(\A\rtimes_r F)$ by using the Green-Julg theorem (see \cite[Theorem 11.1]{GHT2000}, one is also referred to \cite[Theorem 17.A.2]{Gre-Jul} for a lovely proof by S.~Echterhoff). Since $\lim\limits_{R\to\infty}\A_R\rtimes_rF\cong\A\rtimes_rF$, then we have that 
$$\lim\limits_{R\to\infty}K_*(\A_R\rtimes_rF)\cong K_*(\A\rtimes_rF)$$
by using the continuity of $K$-theory. This completes the proof.
\end{proof}

For each $d>0$, we deonte $P_d(\Ga)$ the Rips complex of $\Ga$ at scale $d$. One is referred to Section 7.2 in \cite{HIT2020} for the metric properties of $P_d(\Ga)$. We only mention that the $\Ga$-action on $P_d(\Ga)$ is cocompact, isometric and proper. Let $\EG$ be the classifying space for proper $\Ga$-action. It is proved in \cite[Theorem 7.4.8]{HIT2020} that $\bigcup_{d\to\infty}P_d(X)$ is a model of $\EG$, where we view $P_d(X)$ as a subspace of $P_{d'}(X)$ for $d\leq d'$. Thus, we have the following commuting diagram:
$$\begin{tikzcd}
\lim\limits_{d\to\infty}K_*(C^*_L(P_d(\Ga),\A)^{\Ga}) \arrow[r, "ev_*"] \arrow[d, "\cong "'] &\lim\limits_{d\to\infty}K_*(C^*(P_d(\Ga),\A)^{\Ga}) \arrow[d, "\cong"] \\
RK^{\Ga}_*(\EG,\A) \arrow[r, "\mu_{\A}"]                       & K_*(\A\rtimes_r\Ga).                   
\end{tikzcd}$$

\section{Proof of the main theorem}

We shall prove this theorem in two steps. Let $F\leq\Ga$ be a finite subgroup, $V\subseteq M$ a $F$-invariant bounded subset. Assume that $\Ga\cdot V$ is homeomorphic to $\Ga\times_FV$. Denote $W=\Ga\cdot V$, $\A_V=\overline{C_0(V)\cdot\A}$ and $\A_W=\overline{C_0(W)\cdot\A}$. Notice that $\A_W$ is a $\Ga$-algebra and $\A_V$ is a $F$-algebra. First of all, we shall prove
$$ev_*:\lim_{d\to\infty}K_*(C^*_L(P_d(\Ga),\A_W)^{\Ga})\to \lim_{d\to\infty}K_*(C^*(P_d(\Ga),\A_W)^{\Ga})$$
is an isomorphism.

\begin{Lem}\label{Roe orbit}
With notations as above, for each $d>0$, we have the following isomorphism:
$$C^*(P_d(\Ga),\A_W)^{\Ga}\cong \A_V\rtimes_rF.$$
\end{Lem}

\begin{proof}
For any $f\in C_b(M)$, then $f$ defines a multiplier on $\K\ox\A$ by $f\cdot(k\ox a)=k\ox fa$. For any $T\in\IC[P_d(\Ga),\A_W]^{\Ga}$, we define $T_V=(1\ox\chi_V)T$, i.e,
$$T_V(x,y)=\chi_V T(x,y),$$
where $\chi_V$ is a bounded continuous function on $W=\Ga\cdot V$. Since $P_d(\Ga)$ is cobounded, there exists a bounded subset $K_0\subseteq P_d(\Ga)$ such that $P_d(\Ga)=\Ga\cdot K_0$. Set
$$K=\{x\in P_d(\Ga)\mid d(x,K_0)\leq\prop(T)\}.$$
Then the set $S=\{(x,y)\in K\times K\mid T(x,y)\ne0\}$ is finite. By Lemma \ref{support condition}, we can assume that $\supp(T(x,y))$ is uniformly bounded. Then $\supp(T(x,y))\cap \gamma V\ne\emptyset$ for only finite many $\gamma\in\Ga$. Then for any $(x,y)\in S$, there are only finite many elements in the orbit of $\Ga\cdot x$ such that $T_V(x,y)\ne 0$. Since $T(x,y)\ne 0$ means that $(x,y)$ must belong to $\Ga\cdot S=\{(\gamma x,\gamma y)\mid(x,y)\in S\}$, we conclude that there are only finite $(x,y)\in Z_d\times Z_d$ such that $T_V(x,y)\ne 0$. Thus $T_V$ can be viewed as an $F$-invariant operator over some $F$-invariant bounded subset $B\subseteq P_d(\Ga)$. This means that $T_V\in C^*(B,\A_V)^{F}$. Notice that the correspondence $T\mapsto T_V$ gives a $C^*$-homomorphism
$$C^*(P_d(\Ga),\A_W)^{\Ga}\to \lim_{\substack{B\subseteq P_d(\Ga)\\\text{$F$-invariant, bounded}}}C^*(B,\A_V)^F.$$

Let $B$ be a $F$-invariant bounded subset of $P_d(\Ga)$. For any $T_B\in\IC[B,\A_V]^F$, we define $T\in\IC[P_d(\Ga),\A_W]^{\Ga}$ by
$$T(x,y)=\sum_{[\gamma F]\in\Ga/F}\gamma^{-1}T(\gamma x,\gamma y).$$
The correspondence $T_B\mapsto T$ rise to a $C^*$-homomorphism:
$$\lim_{\substack{B\subseteq P_d(\Ga)\\\text{$F$-invariant, bounded}}}C^*(B,\A_V)^F\to C^*(P_d(\Ga),\A_W)^{\Ga}$$
It is clear that this map is the inverse of the map as above. This shows that these two algebras are isomorphic to each other. Since $C^*(B,\A_V)$ is isomorphic to $\A_V\rtimes_r F$ by Proposition \ref{cp and Roe} for any bounded $B$, this completes the proof.
\end{proof}

By using a basicly same argument, one can also prove the analogous statement for localization algebras:

\begin{Lem}\label{local orbit}
Let $F$ be a finite subgroup of $\Ga$. Then
$$\lim_{d\to\infty}K_*(C^*_L(P_d(\Ga),\A_W)^{\Ga})\cong K_*(C^*_L(P_d(F),\A_V)^F).$$
\end{Lem}

\begin{proof}
By using an argument as Lemma \ref{Roe orbit}, one can also prove that
$$\lim_{\substack{B\subseteq P_d(\Ga)\\\text{$F$-invariant, bounded}}}C^*_L(B,\A_V)^F\cong C^*_L(P_d(\Ga),\A_W)^{\Ga}.$$
After taking the inductive limit on $d$, the left side contains all $P_r(F)$ for $r\geq 0$. Thus, the left side is also made into a classifying space for proper $F$-action. We then conclude the lemma.
\end{proof}

\begin{Pro}\label{iso orbit}
With notations as above, the evaluation map
$$ev_*:\lim_{d\to\infty}K_*(C^*_L(P_d(\Ga),\A_W)^{\Ga})\to \lim_{d\to\infty}K_*(C^*(P_d(\Ga),\A_W)^{\Ga})$$
is an isomorphism.
\end{Pro}

\begin{proof}
By Lemma \ref{Roe orbit} and Lemma \ref{local orbit}, we have the following commuting diagram:
$$\begin{tikzcd}
\lim\limits_{d\to\infty}K_*(C^*_L(P_d(\Ga),\A_W)^{\Ga}) \arrow[r, "ev_*"] \arrow[d, "\cong "'] &\lim\limits_{d\to\infty}K_*(C^*(P_d(\Ga),\A_W)^{\Ga}) \arrow[d, "\cong"] \\
\lim\limits_{d\to\infty}K_*(C^*_L(P_d(F),\A_V)^{F}) \arrow[r, "ev_*"]                       & \lim\limits_{d\to\infty}K_*(C^*(P_d(F),\A_V)^{F}).                   
\end{tikzcd}$$
Notice that the classifying space for proper $F$-action can be chosen to a single point since $F$ is finite. Then $\lim\limits_{d\to\infty}K_*(C^*_L(P_d(F),\A_V)^{F})$ is isomorphic to $K_*(C^*_L(pt,\A_V)^F)$, which is isomorphic to $K_*(C^*(pt,\A_V)^F)$ by evaluation map (see \cite[Lemma 12.4.3]{HIT2020}). The right bottem term of the diagram above is clearly isomorphic to $K_*(\A_V\rtimes_rF)$. This completes the proof.
\end{proof}

Let $W_i=\Ga\cdot V_1$ and $W_2=\Ga\cdot V_2$ for two bounded sets $V_1,V_2\subseteq M$. Define $\A_{W_1}$ and $\A_{W_2}$ as above. The following lemma is similar with \cite[Lemma 6.3]{Yu2000}.

\begin{Lem}\label{cutting coeff}
Assume that $W=W_1\cup W_2$ and $W'=W_1\cap W_2$, then for each $d\geq 0$,\begin{itemize}
\item[(1)] $K_*(C^*(P_d(\Ga),\A_{W_1})^{\Ga})+K_*(C^*(P_d(\Ga),\A_{W_2})^{\Ga})=K_*(C^*(P_d(\Ga),\A_W)^{\Ga})$;
\item[(2)] $K_*(C^*(P_d(\Ga),\A_{W_1})^{\Ga})\cap K_*(C^*(P_d(\Ga),\A_{W_2})^{\Ga})=K_*(C^*(P_d(\Ga),\A_{W'})^{\Ga})$;
\item[(3)] $K_*(C^*_L(P_d(\Ga),\A_{W_1})^{\Ga})+K_*(C^*_L(P_d(\Ga),\A_{W_2})^{\Ga})=K_*(C^*_L(P_d(\Ga),\A_W)^{\Ga})$;
\item[(4)] $K_*(C^*_L(P_d(\Ga),\A_{W_1})^{\Ga})\cap K_*(C^*_L(P_d(\Ga),\A_{W_2})^{\Ga})=K_*(C^*_L(P_d(\Ga),\A_{W'})^{\Ga})$.
\end{itemize}\end{Lem}

\begin{proof}
The proof is an equivariant version of \cite[Lemma 6.3]{Yu2000}. One is also referred to \cite[Lemma 3.5]{SW2007}. One can choose a equivariant partition of unity $\{\phi_1,\phi_2\}$ on $W$ such that $\supp(\phi_i)\subseteq W_i$ for $i=1,2$ instead of $\varphi_i$ in \cite[Lemma 3.5]{SW2007}.
\end{proof}

\begin{proof}[Proof of Theorem \ref{main result}]
For each $R>0$, denote $B_0(R)=B(m_0,R)\subseteq M$, where $m_0$ is the based point of $M$. Denote $B(R)=\Ga\cdot B_0(R)$ and $\A_{B(R)}=\overline{C_0(B(R))\cdot\A}$. It is clear that
$$\lim_{R\to\infty}C^*_L(P_d(\Ga),\A_{B(R)})^{\Ga}=C^*_L(P_d(\Ga),\A)^{\Ga}\qquad \lim_{R\to\infty}C^*(P_d(\Ga),\A_{B(R)})^{\Ga}=C^*(P_d(\Ga),\A)^{\Ga}.$$
Then we have the following commuting diagram:
$$\begin{tikzcd}
\lim\limits_{d\to\infty}K_*(C^*_L(P_d(\Ga),\A)^{\Ga}) \arrow[r] \arrow[d, "\cong "'] & \lim\limits_{d\to\infty}K_*(C^*_L(P_d(\Ga),\A)^{\Ga}) \arrow[d, "\cong"] \\
\lim\limits_{d\to\infty}\lim\limits_{R\to\infty}K_*(C^*_L(P_d(\Ga),\A_{B(R)})^{\Ga}) \arrow[r] \arrow[d, "\cong "']  & \lim\limits_{d\to\infty}\lim\limits_{R\to\infty}K_*(C^*_L(P_d(\Ga),\A_{B(R)})^{\Ga})   \arrow[d, "\cong"]                 \\
\lim\limits_{R\to\infty}\lim\limits_{d\to\infty}K_*(C^*_L(P_d(\Ga),\A_{B(R)})^{\Ga})\arrow[r]     &    \lim\limits_{R\to\infty}\lim\limits_{d\to\infty}K_*(C^*_L(P_d(\Ga),\A_{B(R)})^{\Ga}),
\end{tikzcd}$$
the order of these two limits can be exchanged since
$$\begin{tikzcd}
K_*(AL(d,R)) \arrow[rr] \arrow[rd] \arrow[dd] && K_*(A(d,R)) \arrow[rd] \arrow[dd] &              \\
 & K_*(AL(d',R)) \arrow[rr] \arrow[dd] &&K_*(A(d',R)) \arrow[dd] \\
K_*(AL(d,R')) \arrow[rd] \arrow[rr]   && K_*(A(d,R')) \arrow[rd]            &              \\
 & K_*(AL(d',R')) \arrow[rr]  && K_*(A(d',R'))           
\end{tikzcd}$$
where $AL(d,R)=C^*_L(P_d(\Ga),\A_{B(R)})^{\Ga}$ and $A(d,R)=C^*(P_d(\Ga),\A_{B(R)})^{\Ga}$ for $d'>d\geq 0$ and $R'>R>0$.

Then it suffices to prove
$$ev_*:\lim_{d\to\infty}C^*_L(P_d(\Ga),\A_{B(R)})^{\Ga}\to \lim_{d\to\infty}C^*(P_d(\Ga),\A_{B(R)})^{\Ga}$$
is an isomorphism for each $R>0$. Since $B(R)$ is cobounded, which is also cocompact. Then there exists a finite cover $\{W_1,\cdots,W_N\}$ of $B(R)$ such that each $W_i$ is of the form $\Ga\times_{\Ga_i}V_i$ for some finite subgroup $\Ga_i\leq\Ga$. By Lemma \ref{cutting coeff} and five lemma, it suffices to prove
$$ev_*:\lim_{d\to\infty}C^*_L(P_d(\Ga),\A_{W_i})^{\Ga}\to\lim_{d\to\infty} C^*(P_d(\Ga),\A_{W_i})^{\Ga}$$
is an isomorphism which has been proved in Proposition \ref{iso orbit}. We then finish the proof.
\end{proof}

\bibliographystyle{plain}
\bibliography{ref}

\begin{thebibliography}{10}

\bibitem{Black}
B.~Blackadar.
\newblock {\em {$K$}-theory for operator algebras}, volume~5 of {\em
  Mathematical Sciences Research Institute Publications}.
\newblock Cambridge University Press, Cambridge, second edition, 1998.

\bibitem{GWY2008}
G.~Gong, Q.~Wang, and G.~Yu.
\newblock Geometrization of the strong {N}ovikov conjecture for residually
  finite groups.
\newblock {\em J. Reine Angew. Math.}, 621:159--189, 2008.

\bibitem{GHT2000}
E.~Guentner, N.~Higson, and J.~Trout.
\newblock Equivariant {$E$}-theory for {$C^*$}-algebras.
\newblock {\em Mem. Amer. Math. Soc.}, 148(703):viii+86, 2000.

\bibitem{HK2001}
N.~Higson and G.~Kasparov.
\newblock {$E$}-theory and {$KK$}-theory for groups which act properly and
  isometrically on {H}ilbert space.
\newblock {\em Invent. Math.}, 144(1):23--74, 2001.

\bibitem{Gre-Jul}
D.~Husem\"{o}ller, M.~Joachim, B.~Jur\v{c}o, and M.~Schottenloher.
\newblock {\em Basic bundle theory and {$K$}-cohomology invariants}, volume 726
  of {\em Lecture Notes in Physics}.
\newblock Springer, Berlin, 2008.
\newblock With contributions by Siegfried Echterhoff, Stefan Fredenhagen and
  Bernhard Kr\"{o}tz.

\bibitem{KS2003}
G.~Kasparov and G.~Skandalis.
\newblock Groups acting properly on ``bolic'' spaces and the {N}ovikov
  conjecture.
\newblock {\em Ann. of Math. (2)}, 158(1):165--206, 2003.

\bibitem{KY2012}
G.~Kasparov and G.~Yu.
\newblock The {N}ovikov conjecture and geometry of {B}anach spaces.
\newblock {\em Geom. Topol.}, 16(3):1859--1880, 2012.

\bibitem{Kas88}
G.~G. Kasparov.
\newblock Equivariant {$KK$}-theory and the {N}ovikov conjecture.
\newblock {\em Invent. Math.}, 91(1):147--201, 1988.

\bibitem{SW2007}
L.~Shan and Q.~Wang.
\newblock The coarse geometric {N}ovikov conjecture for subspaces of
  non-positively curved manifolds.
\newblock {\em J. Funct. Anal.}, 248(2):448--471, 2007.

\bibitem{Tu1999}
J.-L. Tu.
\newblock The {B}aum-{C}onnes conjecture and discrete group actions on trees.
\newblock {\em $K$-Theory}, 17(4):303--318, 1999.

\bibitem{HIT2020}
R.~Willett and G.~Yu.
\newblock {\em Higher index theory}.
\newblock Cambridge University Press, 2020.

\bibitem{Yu1997}
G.~Yu.
\newblock Localization algebras and the coarse {B}aum-{C}onnes conjecture.
\newblock {\em $K$-Theory}, 11(4):307--318, 1997.

\bibitem{Yu2000}
G.~Yu.
\newblock The coarse {B}aum-{C}onnes conjecture for spaces which admit a
  uniform embedding into {H}ilbert space.
\newblock {\em Invent. Math.}, 139(1):201--240, 2000.

\end{thebibliography}

\vskip 5cm
\begin{itemize}
\item[] Liang Guo\\
Research Center for Operator Algebras, School of Mathematical Sciences, East China Normal University, Shanghai, 200241, P. R. China.\quad
E-mail: 52205500015@stu.ecnu.edu.cn

\item[] Qin Wang \\
Research Center for Operator Algebras,  and Shanghai Key Laboratory of Pure Mathematics and Mathematical Practice, School of Mathematical Sciences, East China Normal University, Shanghai, 200241, P. R. China. \quad
E-mail: qwang@math.ecnu.edu.cn
\end{itemize}
\end{document}